\input amstex
\input epsfx.tex
\frenchspacing
\documentstyle{amsppt}
\magnification=\magstep1
\baselineskip=14pt
\vsize=18.5cm
\footline{\hfill\sevenrm version 20051125} 
\def\cal{\Cal} 
\def\Tr{{\text{\rm Tr}}}

\def\End{{\text{\rm End}}}

\def\mod{\bmod}
\def\congr{\equiv}

\def\O{{\Cal O}}
\def\Que{\bold Q}
\def\CC{\bold C}
\def\FF{\bold F}
\def\HH{\bold H}
\def\Zee{\bold Z}
\def\Fp{\FF_p}
\def\gothp{{\goth p}}
\newcount\refCount
\def\newref#1 {\advance\refCount by 1
\expandafter\edef\csname#1\endcsname{\the\refCount}}
\newref AKS  %AKS
\newref BHP  %Baker, Harman, Pintz, The difference between consecutive primes II
\newref Bern %Bernstein, prob. AKS
\newref BR   %thesis Reinier
\newref BS   %B&S ANTS VI paper
\newref BW   %Buhler-Wagon
\newref Coh  %Cohen - course in algorithmic number theory
\newref CH   %Couveignes-Henocq
\newref EN   %Enge
\newref MCA  %Moderen Computer Algebra, Gathen, Gerhard
\newref GS   %Gee-Stevenhagen, ANTS IV
\newref Iv   %Ivic zeta-function
\newref KE   %Kedlaya
\newref KSZ  %Konstantinou, Stamatiou, Zaroliagis - Indocrypt 2003
\newref LZ   %Lay-Zimmer
\newref LP   %improvements by Lenstra, Pomerance to AKS
\newref Mo   %Morain
\newref SA   %Satoh
\newref SSK  %Savas, Schmidt, Koc
\newref Sc   %schoof - 1985 paper
\newref SC   %schoof - Counting points on elliptic curves
\newref SI   %Silverman ell curves
\newref ST   %Stevenhagen Tokyo volume
\newref CF   %Cassels and Frohlich
\topmatter
\title 
Constructing elliptic curves in almost polynomial time
\endtitle
\author Reinier Br\"oker, Peter Stevenhagen\endauthor
\address Mathematisch Instituut,
Universiteit Leiden, Postbus 9512, 2300 RA Leiden, The Netherlands.
\email reinier, psh\@math.leidenuniv.nl\endemail
\endaddress
\abstract
We present an algorithm that, on input of an integer
$N\ge1$ together with its prime factorization,
constructs a finite field $\FF$ and an elliptic curve
$E$ over $\FF$ for which $E({\FF})$ has order $N$.
Although it is unproved that this can be done for all $N$, a heuristic 
analysis shows that the algorithm has an expected run time that is polynomial in
$2^{\omega(N)}\log N$, where $\omega(N)$ is the number
of distinct prime factors of $N$.
In the cryptographically relevant case where $N$ is prime, an
expected run time $O((\log N)^{4+\varepsilon})$ can be achieved.
We illustrate the efficiency of the algorithm by constructing
elliptic curves with point groups of order $N=10^{2004}$ and
$N=\text{nextprime}(10^{2004})=10^{2004}+4863$.
\endabstract
\endtopmatter

\document

\head 1. Introduction
\endhead

\noindent
For an elliptic curve $E$ defined over the finite field $\FF_q$ of $q$
elements, the order $N=\#E(\FF_q)$ of the group
of $\FF_q$-rational points of $E$ is an integer in the Hasse interval
$$
{\cal H}_q=[(\sqrt q-1)^2, (\sqrt q+1)^2]=[q+1-2\sqrt q, q+1+2\sqrt q]
\eqno(1.1)
$$
around $q$.
Various point counting algorithms [\Sc, \SA, \KE] have been
developed over the last 20 years that compute $N$ in 
polynomial time from the standard representation of $E$ by a
Weierstrass equation over $\FF_q$.
A natural `inverse problem' to the point counting problem is the
following.
\proclaim
{Problem 1}
Given a finite field $\FF_q$ and an integer $N\in {\cal H}_q$,
find an elliptic curve $E/\FF_q$ for which $E(\FF_q)$ has
order $N$.
\endproclaim\noindent
If $q=p$ is a prime number, than all integers $N\in {\cal H}_p$
arise as the order of an elliptic curve over $\FF_p$, and
a solution to Problem 1 always exists.
For prime powers $q=p^k$ this is not generally true:
the values $N\in{\cal H}_q$ having $N\congr 1\mod p$ can only be
realized by supersingular elliptic curves over $\FF_q$, and these are
in most cases too rare [\SI, Theorem V.4.1] to account for all values 
$N\congr 1\mod p$ in ${\cal H}_q$.
On the other hand, all values $N\not\equiv 1\mod p$ in ${\cal H}_q$
do arise as orders of elliptic curves over $\FF_q$.

No algorithm is known to solve problem 1 (in the cases where a
solution exists) in a time that is polynomially bounded
in the input size $\log q\approx \log N$.
Due to the fact that point counting of elliptic curves over $\FF_q$
can be done in polynomial time, the naive probabilistic algorithm
of trying random curves $E/\FF_q$ until a curve with the right
number of points is found has expected run time $\widetilde O(N^{1/2})$.
Here we use the $\widetilde O$-notation to indicate that terms that are
of logarithmic size in the main term have been disregarded.

Simple-minded as it is, the naive algorithm compares favorably to
the deterministic complex multiplication algorithm to solve
Problem 1 that is discussed in the next section.
This is due to the size of the auxiliary polynomials (`class polynomials')
in that algorithm, which become prohibitively large for most pairs $(q, N)$.
In order to obtain algorithms that are substantially better than
the naive method, one can relax the conditions in Problem 1 in the 
following way.
\proclaim
{Problem 2}
Given an integer $N\ge1$, find a finite field\/ $\FF$ and an elliptic
curve $E/\FF$ for which $E(\FF)$ has order $N$.
\endproclaim\noindent
In the case where the discrete logarithm problem in $E(\FF)$ is the basis
of a cryptosystem, it is important that $N$ has certain properties, e.g.,
that it is divisible by or equal to a large prime number, whereas the precise
value of $q=\#\FF$ is less relevant.
In this case one needs a solution to Problem 2, not to Problem 1.
The observation is not new, and both problems occur in the list
of problems in the introduction of [\LZ] that `can be solved'.

The main result of this paper is that, even though no efficient
solution to Problem~1 is known, Problem 2 does admit such a solution
if $N$ is provided to the algorithm in factored form.
For practical applications, such as those in elliptic curve cryptography,
it is unlikely that one will need or want to use elliptic
curves for which the factorization of the group order is unknown,
so requiring the factorization of~$N$ to be part of the input
is not a severe restriction.
Our solution to Problem~2 for factored orders $N$ is almost polynomial time, 
provided that one is willing to assume a number
of `standard heuristic assumptions' that we will make
explicit in Section~4.
\proclaim
{Main Theorem}
There exists an algorithm that, on input of an integer $N\ge1$ 
together with its factorization, returns a prime number $p$
and an elliptic curve $E/\FF_p$ with $\#E(\FF_p)=N$ whenever such a pair
$(E, p)$ exists.
Under standard heuristic assumptions, a pair $(E, p)$ exists for all
$N$, and the expected run time of the algorithm is polynomial in
$2^{\omega(N)}\log N$.
Here $\omega(N)$ denotes the number of distinct prime factors of $N$.
\endproclaim\noindent
Although the run time in the Main Theorem is not polynomial in the usual
sense, it is polynomial in $\log N$ outside a zero density subset of 
$\Zee_{\geq 1}$ consisting of very smooth input values $N$.
Note that such $N$ are not used in cryptographic applications,
as the discrete logarithm problem in groups of smooth order tends to be easy.
\proclaim
{Corollary}
If the input values $N$ in the Main Theorem are restricted to be prime 
numbers or,
more generally, to be in the density $1$ subset of $\Zee_{\geq 1}$ consisting of
those $N$ having $\omega(N)<2\log\log N$,
then the expected run time is polynomial in $\log N$.
\endproclaim\noindent
The factorization of $N$ is used by the algorithm in the Main Theorem
to reduce square root extractions of small integers modulo~$N$ to square root 
extractions modulo the prime factors of $N$.
It is here that the approximate number $2^{\omega(N)}$ of such roots
enters the run time of the algorithm.
The precise exponents in the run time depend on one's willingness
to accept fast multiplication techniques and
probabilistic subroutines in the algorithm.
For instance, the square root extractions of small integers modulo the
prime factors of $N$ can be done efficiently by probabilistic means or,
much less efficiently, but still in time polynomial in 
$2^{\omega(N)}\log N$, by a deterministic algorithm [\Sc].
Similarly, one may require for the prime number $p$ returned by
the algorithm that its primality is proved by a deterministic
AKS-type polynomial time algorithm, or employ a faster probabilistic
algorithm to do so.
If we insist on guaranteed correct output, i.e., a {\it proven\/} prime
$p$ as the characteristic of our curve~$E$, but allow fast
multiplication and probabilistic subroutines of the kind mentioned above,
the heuristic run time of our algorithm is 
$O(2^{\omega(N)}(\log N)^{4+\varepsilon})$ for every $\varepsilon>0$
(Corollary 4.4.).
In the cryptographically relevant case where $N$ is prime [\SSK, \KSZ],
this becomes $O((\log N)^{4+\varepsilon})$
(Corollary 4.2).

It should not come as a surprise that our solutions to Problem 2
are elliptic curves defined over prime fields.
Indeed, it is easy to see that the union of the Hasse intervals
${\cal H}_q$ over the prime powers $q$ that are {\it not\/} primes
is a zero density subset of~$\Zee_{\ge1}$.
Solvability of Problem 2 for all values of $N$ is therefore 
in an informal sense `equivalent' to the fact that the union
of the Hasse intervals ${\cal H}_p$
over the primes $p$ contains $\Zee_{\ge1}$.
Defining the Hasse interval around arbitrary integers $q$ by formula
(1.1), we have the equivalence
$$
N\in {\cal H}_q \Longleftrightarrow  q\in {\cal H}_N,
\eqno(1.2)
$$
and we see that we want every Hasse interval ${\cal H}_N$ around an
integer $N$ to contain a prime number $p$.
This amounts to the statement that the size of the `gaps' between 
consecutive primes around $N$ does not exceed $4\sqrt N$.
Although prime gaps of this size are not believed to exist,
the best proven upper bound on their size [\BHP] is currently
$O(N^\alpha)$, with $\alpha = .525 >{1\over2}$.
Even under assumption of the generalized Riemann hypothesis, the best
result [\Iv, Theorem 12.10] is only $O(N^{1/2}\log N)$.
This means that we have no proof that Problem 2 is solvable for all $N$,
and already for this reason a rigorous run time analysis for our Main Theorem 
is out of reach.

By the prime number theorem, we expect one out of
every $\log N$ integers around $N$ to be prime, so 
the Hasse interval ${\cal H}_N$ of length $4\sqrt N$ around $N$
will normally contain {\it many\/} primes $p$.
In practice, there is always an abundance of
primes $p$ for which there exist elliptic curves $E/\FF_p$ of order $N$,
and it seems extremely unlikely that the number of primes in ${\cal H}_N$,
which grows `on average' as $\widetilde O(N^{1/2})$, will be zero for
some $N$.
The real task of our algorithm is therefore not so much to find a prime
$p\in {\cal H}_N$, but rather to find a prime $p\in {\cal H}_N$ for
which a curve $E/\FF_p$ of order $N$ can be constructed {\it efficiently\/}.
In Section 2, we show how this leads to a new Problem 3, whose efficient
solution yields an efficient solution of Problem 2.

Section 3 describes an Algorithm that solves our Problem 3
and finds a suitable prime $p\in{\cal H}_N$.
Its heuristic run time is derived in Section 4.
It is based on various unproved but reasonable statements, such as the fact
that random integers in ${\cal H}_N$ will be prime with probability $1/\log N$.
We also present numerical evidence for such unproved statements.
In the case where $N$ is prime,
the heuristic arguments are very similar to those going into the
analysis of the elliptic curve primality proving algorithm ECPP [\Mo].

Section 5 comments on an efficient implementation of the
Algorithm to solve Problem 2.
It illustrates its practical applicability 
by treating as examples `random' values of $N$
such as $N=10^{2004}$ and $N=\text{nextprime}(10^{2004})=10^{2004}+4863$.

\head 2. Complex multiplication constructions
\endhead

\noindent
Although much in this section generalizes to arbitrary prime
powers $q$, we now focus on the case relevant to us, where
$q=p>3$ is a prime number and $N\in {\cal H}_p$ an integer that
we want to realize as the order of some elliptic curve $E/\FF_p$.

Constructing an elliptic curve $E/\FF_p$ having $N$ points
roughly comes down to computing the $j$-invariant $j(E)\in\FF_p$
of such a curve, and the theory of complex multiplication
provides a deterministic way of doing so.
If we write $N=p+1-t$, then $E/\FF_p$ has $\#E(\FF_p)=N$ if and only if
the Frobenius morphism $F_p$ of $E$ satisfies the 
quadratic relation 
$$
F_p^2-tF_p+p=0\eqno(2.1)
$$
of discriminant $\Delta = t^2-4p<0$ in $\End(E)$. 
If $F_p$ satisfies (2.1),
then $\Zee[F_p]\subset \End(E)$ is isomorphic to
the imaginary quadratic order $\O_\Delta$ of discriminant $\Delta$,
and $F_p$ corresponds to the element
$\pi={t+\sqrt\Delta\over2}\in \O_\Delta$ of trace $t$ and norm $p$.
Unless we are in the supersingular case $t=0$ having $\Delta=-4p$,
which is too special to be of interest here,
this means that
$p=\pi\bar\pi$
splits into principal primes in $\O_\Delta$.

Over the the field $\CC$ of complex numbers, it is a classical result that 
the isomorphism classes of elliptic curves having endomorphism ring 
isomorphic to $\O_\Delta$ correspond to the classes of invertible 
$\O_\Delta$-ideals in the class group $\text{Pic}(\O_\Delta)$ 
of the order $\O_\Delta$.
Invertible $\O_\Delta$-ideals can be viewed as lattices in $\CC$, and
the $j$-invariants of these lattices are precisely the $j$-invariants
of the elliptic curves having endomorphism ring isomorphic to $\O_\Delta$.
It follows that we can evaluate these $j$-invariants as values of the
modular function $j: \HH \rightarrow \CC$ in points $\tau_Q$ in the
complex upper half plane $\HH$ representing the ideal classes
$[Q]\in\text{Pic}(\O_\Delta)$.
More precisely,
if we represent the ideal classes of $\text{Pic}(\O_\Delta)$ in the
standard way [\Coh, Section 5.2] as reduced binary quadratic forms 
$Q=aX^2+bXY+cY^2$ of discriminant $b^2-4ac=\Delta$, we have
$\tau_Q={-b+\sqrt\Delta\over 2a}$.
The {\it class polynomial\/}
$$
P_\Delta=\prod_{[Q]\in\text{Pic}(\O_\Delta)} (X-j(\tau_Q))\in\Zee[X]
$$
has integer coefficients, so it can be computed exactly from
complex approximations of the $j(\tau_Q)$.
In the ordinary case $t\ne0$, the reduction modulo $p$ of
the class polynomial $P_\Delta$ splits into $h(\Delta)=\#\text{Pic}(\O_\Delta)$
distinct linear factors in $\FF_p[X]$, and the roots are the $j$-invariants 
of the elliptic curves
over $\FF_p$ having endomorphism ring isomorphic to $\O_\Delta$.
If $j_0\ne 0, 1728$ is one of these zeroes in $\FF_p$, then the
curve $E_a/\FF_p$ with Weierstrass equation $Y^2=X^3+aX-a$
has $j$-invariant $j_0$ if we choose $a$ to satisfy
$$
j_0=1728 \,{4a\over 4a+27},
$$
and its number of points is either $N=p+1-t$ or $p+1+t$.
We easily check in which case we are, by point counting or by simply
evaluating $N\cdot P$ and $(p+1+t)\cdot P$ for the point $P=(1,1)$ on $E_a$.
If the order is $N$ we are done;
if not, then the quadratic twist $Y^2=X^3+ag^2X-ag^3$ of $E_a$,
with $g$ a non-square in $\FF_p^*$, solves our problem.
In the special cases $j_0=0, 1728$ that we disregard here,
there are a few more quadratic twists to consider --
see Example 5.2.

Most of the work in the complex multiplication method
goes into the computation of the class polynomial $P_\Delta$.
As the degree of $P_\Delta$ and the size of its coefficients
both grow like $|\Delta|^{1/2}$ for $\Delta\to-\infty$,
the run time can be no better than $\widetilde O(|\Delta|)$.
This is the actual run time [\EN] for the classical analytic approach
using the modular function $j: \HH \rightarrow \CC$.
The same is true for the more recent non-archimedean approach [\BS, \CH]
to the evaluation of $P_\Delta$, which approximates the roots of 
$P_\Delta$ by a Newton iteration process over $\Que_\ell$
for a suitable small prime $\ell$.
For both methods, it is possible to reduce the run time by
sizable constant factors if one replaces the $j$-function by
`smaller' modular functions [\GS, \ST, \BR].
This is very important from a practical, but not from a computational
complexity point of view.

In the complex multiplication method, one can save some work by
computing the class polynomial $P_D$ for the fundamental discriminant
$D=\text{disc}(\Que(\sqrt\Delta))$ rather than that for $\Delta$ itself.
As $p=\pi\bar\pi\in \O_\Delta$ splits in the same way in the maximal order
$\O_D\supset \O_\Delta$ as it does in $\O_\Delta$,
elliptic curves over $\FF_p$ with endomorphism ring
$\O_D$ are just as good for our purposes, and we may everywhere replace
$\Delta$ by $D$ in the algorithm.
If $\Delta$ has a large square factor, this can be a considerable improvement
since the polynomial $P_D$ is then much smaller than $P_\Delta$.

If we apply the complex multiplication method to solve Problem 1, we have
no control over the discriminant 
$$
\Delta=\Delta(p, N)=t^2-4p=(p+1-N)^2-4p, \eqno(2.2)
$$
which will typically be of the same order of magnitude as $N$ and
without large square factors.
In that case, the resulting run time $\widetilde O(N)$ is inferior to
the $\widetilde O(N^{1/2})$ of the naive probabilistic method.

For Problem 2, the situation is different as only $N$ is then
given as input, and we typically have {\it many\/} primes
$p\in {\cal H}_N$ to choose from.
An obvious thing to do here is to choose $p\in {\cal H}_N$ as close
as possible to the end points of the interval, so that the absolute value 
of the trace $t=p+1-N$ differs from $2\sqrt p$ by a small amount.
By the prime number theorem, we expect to be able to find $p$ for which
$|t|-2\sqrt p$ is of size $\log N$.
This makes $\Delta=t^2-4p$ of size $\widetilde O(N^{1/2})$, and reduces the
run time of the algorithm to $\widetilde O(N^{1/2})$, just as for the
naive probabilistic method.

More generally, one can examine which primes $p$ at distance at most
$N^\alpha$ from the end points of ${\cal H}_N$ give rise to values
of $\Delta$ with large square factors.
Heuristically, there are about $N^\alpha/\log N$ such primes, giving
rise to discriminants of size $N^{\alpha+1/2}$.
Among the discriminants of this size, those of the form $\Delta=f^2D$ with
$|D|<N^\beta$ constitute a fraction of order of magnitude
$$
P(\alpha, \beta)
  =N^{-(\alpha+1/2)}
   \sum_{|D|<N^\beta \text{ squarefree}} \sqrt{N^{\alpha+1/2}\over |D|}
       \approx N^{{1\over 2}(\beta-\alpha)-{1\over 4}}.$$
The number of discriminants $\Delta=f^2 D$ with $|D|<N^\beta$
we expect to find from $p$'s no further than $N^\alpha$ from the end 
points of ${\cal H}_N$ is therefore 
$$
P(\alpha, \beta)\cdot{N^\alpha\over\log N}=
{1\over\log N}\cdot
N^{{1\over 2}(\alpha+\beta)-{1\over 4}}
,$$
which tends to infinity with $N$ exactly when we have 
$\alpha+\beta > 1/2$.
Rough as this heuristic analysis may be, it `explains' why in the
example $N=10^{30}$ given in [\BS, Section 6] to illustrate the
non-archimedean approach to computing class polynomials, examining the primes
$p$ at distance $<10^6$ from the end points of ${\cal H}_N$ leads to a
fundamental discriminant $D\approx -10^8$.
As examining the primes in an interval of length $N^\alpha$ to achieve
$|D|<N^\beta$ gives rise to a run time 
$\widetilde O(N^{\max\{\alpha, \beta\}})$,
we can achieve a heuristic run time $O(N^{{1\over4}+\varepsilon})$ by
taking $\alpha=\beta={1\over 4}+\varepsilon$.
Although this is still exponential, this method of selecting $p$
already enables us to deal with values of $N$ the naive method cannot handle.

The extreme case $(\alpha, \beta)=(\varepsilon, 1/2)$ corresponds to taking $p$ as
close as possible to the end points of ${\cal H}_N$, a case we 
already discussed.
The other extreme $(\alpha, \beta)=(1/2, \varepsilon)$ indicates that it should
be possible to find $D$ of {\it subexponential\/} size in terms of our
input length $\log N$. 
This suggests that a fruitful approach to solving
Problem 2 by the complex multiplication method consists in efficiently
minimizing the fundamental discriminant $D$ involved.

It turns out that we can actually determine
the `minimal' imaginary quadratic fundamental discriminant $D$ that can be
used to construct an elliptic curve of order $N$ in a relatively
straightforward way.
It uses the `symmetry' between the order $N$ of the
point group $E(\FF)$ and the order $q=p$ of $\FF$ itself, which are norms
of the quadratic integers $1-\pi=1-F_p$ and $\pi=F_p$, respectively.
This symmetry is already familiar to us from (1.2).
In the case of the discriminant 
$\Delta=(\pi-\bar\pi)^2=((1-\pi)-(1-\bar\pi))^2$ in (2.2),
it takes the form
$$
\Delta(p, N)=(p+1-N)^2-4p=(N+1-p)^2-4N.
$$
We now fix $N$ and try to write $\Delta=\Delta(p)$ as
$$
\Delta(p)=(N+1-p)^2-4N= f^2D \eqno(2.3)
$$
for `small' $D<0$.
This comes down to solving the positive definite equation
$$
x^2-Df^2 = 4N\eqno(2.4)
$$
in integers $x$ and $f$ in such a way that the number
$p=N+1-x$ is prime.
This leads us to the following problem.
\proclaim
{Problem 3}
Given an integer $N\ge1$, find the smallest squarefree integer $d\ge1$ 
together with an algebraic integer $\alpha\in K=\Que(\sqrt{-d})$ such
that
\itemitem{\rm (i)} 
$N_{K/\Que}(\alpha)=N$;
\itemitem{\rm(ii)} 
$p=N_{K/\Que}(1-\alpha)=N+1-\Tr_{K/\Que}(\alpha)$ is prime.
\endproclaim\noindent
The prime $p$ occurring in condition (ii) has the property that
there exists an elliptic curve $E/\FF_p$ having $N$ points
and endomorphism ring $\End(E)$ isomorphic to the ring of integers 
$\O_K$ of $K=\Que(\sqrt{-d})$.
Once we find the solution $(\alpha, d)$ to Problem~3,
we can use it to solve Problem 2 for that same $N$:
take $p=N_{K/\Que}(1-\alpha)$ and construct an elliptic curve over $\FF_p$ with
endomorphism ring $\O_K$ for which $1-\alpha\in\O_K$ is the Frobenius,
using the class polynomial for the order $\O_K$.
This elliptic curve will have $N=N_{K/\Que}(\alpha)$ points, as desired.
\head 3. An Algorithm to solve Problem 3
\endhead

\noindent
As indicated in the introduction, it is not possible to prove
rigorously that {\it any\/} pair $(\alpha, d)$ meeting the conditions
of  Problem 3 exists at all,
let alone that there is a pair with small $d$ that can be found efficiently.
We will however argue in the next section why it is reasonable to expect
that the smallest integer $d$ solving Problem~3 exists for all $N\ge1$, and
why this $d$ is even rather small in terms of $N$, of size at most
$\widetilde O((\log N)^2+2^{\omega(N)})$.
Given this expectation, it makes sense to solve Problem 3 in a
straightforward way using an algorithm that,
on input of a factored number $N$, tries for increasing squarefree numbers
$d\in\Zee_{\geq 1}$ to
\item{--}
find the integral ideals in $K=\Que(\sqrt{-d})$ of norm $N$;
\item{--}
determine the generators of those ideals that are principal;
\item{--}
test for each generator $\alpha$ found whether $N_{K/\Que}(1-\alpha)$ is prime.
\vskip0cm\noindent
As soon as a prime value $p=N_{K/\Que}(1-\alpha)$ is encountered
for some $d$, this is the minimal $d$ we are after, and $(\alpha, d)$ is
a solution to Problem 3.

Before we describe an actual algorithm, we look at the three
individual tasks to be performed, and the run time of the various
subroutines involved.
These run times depend on the time $O(L^{1+\mu})$ needed to multiply
two $L$-bit integers.
We have $\mu=1$ for ordinary multiplication, and $\mu=\varepsilon>0$
for any fast multiplication method.
We will give our run times using $\mu=\varepsilon>0$.
\medskip\noindent
{\it Task 1: Finding the integral ideals in
$\Que(\sqrt{-d})$ of norm $N$.}
\smallskip\noindent
Write the ring of integers of $\Que(\sqrt{-d})$ as $\Zee[\omega]$,
with $\omega=\omega_d$ a zero of
$$
f=f^\omega_\Que=
\cases
X^2-X+{1+d\over 4}&\text{if $-d\congr1\mod 4$;}\cr
X^2+d&\text{otherwise.}\cr
\endcases
\eqno{(3.1)}
$$
Then every ideal of norm $N$ in $\Zee[\omega]$ can uniquely be written as
$kI$, with $k$ a positive integer for which $k^2$ divides $N$,
and $I$ a {\it primitive\/} ideal of $\Zee[\omega]$ of norm $N_0=N/k^2$.
This last condition means that $\Zee[\omega]/I$ is cyclic
of order $N_0$, and it implies that we have $I=(N_0, \omega-r)$ for some
integer $r\in\Zee$ satisfying $f(r)\congr0\mod N_0$.
Finding all ideals of norm $N$ therefore amounts to finding, for each
square divisor $k^2|N$, the roots of $f$ modulo $N_0=N/k^2$.
It is here that we need to have the factorization of $N$ at our disposal,
not only because this implicitly encodes a list of square divisors $k^2|N$,
but also because it enables us to find the roots of $f$ modulo~$N_0$.
Indeed, finding these roots is done by finding the roots of $f$ modulo
the prime powers $p^{\text{ord}_p(N_0)}$ dividing $N_0$, and combining
these in all possible ways, using the Chinese remainder theorem, to obtain
the roots modulo $N_0$.
Note that $f$ has {\it no\/} roots modulo $N_0$ if $N_0$ is divisible
by a prime $p$ that is inert in $\Zee[\omega]$, or by the square $p^2$
of a prime $p$ that ramifies in $\Zee[\omega]$.

As finding a root of $f$ modulo an integer essentially amounts to
extracting a square root of $-d$ modulo that integer, we need to extract
square roots of $-d$ modulo the prime powers dividing $N_0$.
This easily reduces to extracting square roots of $-d$ modulo each of the
primes dividing $N_0$.
This can be done efficiently by employing a variant of the
(probabilistic) Cantor-Zassenhaus algorithm [\MCA, Section 14.5],
and leads to an expected run time $O((\log p)^{2+\varepsilon})$ to extract 
square roots modulo a prime $p$.
For any selection of square roots  $(\sqrt{-d} \mod p^{\text{ord}_p(N_0)})$,
the Chinese remainder theorem lifts these to a square root modulo $N_0$
in time $O(\omega(N)(\log N)^2)$.
\medskip\noindent
{\it Task 2: Finding generators for principal ideals of norm $N$.}
\smallskip\noindent
For each ideal $kI=k\cdot(N_0, \omega-r)\subset \Zee[\omega]$
of norm $N$ found, we use the 1908 algorithm of Cornacchia
described in [\SC, pp. 229--232] or [\BW] to find a generator of~$I$,
if it exists.
This algorithm performs a number of steps of the
Euclidean algorithm to the basis elements $N_0$ and $\omega-r$
of the $\Zee$-lattice $I=(N_0, \omega-r)\subset \Zee[\omega]$ in order
to decide whether $I$ is a principal ideal.
If it is, a generator $\alpha=k\alpha_0$ of $kI$ of norm $N$ is found.
The other generator of $I$ is $-\alpha$.
For the special values $d=1$ and $d=3$ there are 4 and 6 generators for each
principal ideal $I$, respectively, obtained by multiplying $\alpha$
by 4th and 6th roots of unity.
The run time of Cornacchia's algorithm on input $k\cdot(N_0, \omega-r)$
is of order $O((\log N )^{2+\varepsilon})$. 
\medskip\noindent
{\it Task 3: Testing which algebraic integers $\alpha$ of norm $N$
lead to prime elements $1-\alpha$.}
\smallskip\noindent
For each of the elements $\alpha$ of norm $N$ found in the previous 
step 2, we need to test whether the norm $N+1-\Tr(\alpha)$
of $1-\alpha$ is a prime number.
As most $\alpha$'s will have norms that are not prime,
a cheap compositeness test such as the Miller-Rabin test
(which takes time $\widetilde O(\log N)$) can be used to discard most $\alpha$'s.
Once we find $\alpha$ for which $N+1-\Tr(\alpha)$ is a probable prime,
we do a true primality test to {\it prove\/} primality of $p=N+1-\Tr(\alpha)$.
This can be done deterministically in time polynomial in $\log N$ by
the 2002 result of  Agrawal, Kayal and Saxena [\AKS].
Recent speed-ups of the test [\LP] take time $O((\log N)^{6+\varepsilon})$,
whereas probabilistic versions [\Bern] have expected run time 
$O((\log N)^{4+\varepsilon})$.
\medskip\noindent
Using the various subroutines specified in the tasks above, we 
formulate an Algorithm to solve Problem 3.
A slightly more practical algorithm that we use to actually find
elliptic curves with a given number of points does not exactly follow
the outline below; it is discussed in Section 5.
The version in this section is phrased to facilitate the heuristic run time
estimate in Section 4.
\medskip\noindent
{\bf Algorithm.}
\hfil\break
Input: a factored integer $N=\prod_{i=1}^t p_j^{e_j}$.
Output: a solution $(d, \alpha)$ to Problem 3.
\vskip5pt\parindent=.6cm
\item{\bf 1.} Put $d \leftarrow 1$.
\item{\bf 2.} If $d$ is not squarefree, put $d \leftarrow d+1$ and go to step 2.
Otherwise, define $\omega = \omega_d$ and $f = f_\Que^\omega$ as in (3.1).
\item{\bf 3.} Determine the splitting behavior in $\Zee[\omega]$ of all
prime divisors of $N$.
\itemitem{\bf 3a.} For every prime divisor $p_i$ of $N$ that is inert in
$\Zee[\omega]$, put
$$
k_1 \leftarrow k_1 p_i^{\lfloor {e_i/2}\rfloor}
$$
in case $e_i$ is even. In case $e_i$ is odd, put
$d \leftarrow d+1$ and go to step 2.
\itemitem{\bf 3b.} For every prime divisor $p_i$ of $N$ that ramifies in
$\Zee[\omega]$, put
$$
k_1 \leftarrow k_1 p_i^{\lfloor {e_i/2}\rfloor}.
$$
\item{\bf 4.} Put $N_1 \leftarrow N/k_1^2$.
For every root $(r\mod N_1)$ of $f$ and for
every square divisor $k_2^2 \mid N_1$ do the following.
\itemitem{\bf 4a.} Put $k \leftarrow k_1 k_2$ and $N_0 \leftarrow N/k^2 =
N_1/k_2^2$.
Use Cornacchia to find a generator of $(N_0,\omega-r)
\subset \Zee[\omega]$, in case it exists.
\itemitem{\bf 4b.} If a generator is found, test for all (2, 4 or 6)
generators $\alpha_0$ whether the norm $N+1-{\Tr}(k\alpha_0)$ of
$k\alpha_0\in\Zee[\omega]$ is prime. If it is, return $d$ and
$\alpha=k\alpha_0$ and halt.
\item{\bf 5.} Put $d \leftarrow d+1$ and go to step 2.
\medskip\noindent
The determination of the splitting behavior of the primes $p_i|N$
in $\Zee[\omega]$ in Step~3 amounts to computing the Kronecker symbol
$D\overwithdelims(){p_i}$ for $D=\text{disc}(\Que(\sqrt{-d}))$.
For $p>2$ this is simply the Legendre symbol, which is easily
evaluated by combining quadratic reciprocity with the Euclidean
algorithm.
The factor $k_1$ computed in this step is the minimal `imprimitivity factor'
dividing all ideals of norm $N$ in $\Zee[\omega]$.
It reflects the fact that primitive ideals are not divisible by
inert primes, or by squares of ramified primes.

The evaluation of the roots of $f$ modulo $N_1$ in Step 4 
is done by evaluating the roots of $f$ modulo the various prime powers
dividing $N_1$, and combining these in all possible ways using the
Chinese remainder theorem.
For the ramified primes $p_i$ dividing $N_1$, which occur with exponent 1,
there is a unique (double) root of $f$ modulo $p$.
For splitting primes $p_i$, the polynomial $f$ has exactly 2 different roots
modulo $p_i$, and these lift uniquely to $\Zee_p$.
Finding the roots of $f$ modulo these $p_i$ is non-trivial as
it involves the extraction of a square root $\sqrt{-d}$ modulo $p_i$.
Refining these roots to roots modulo $p_i^{e_i}$ is much faster, and
an easy application of Hensel's lemma.
The number of distinct roots modulo $N_1$ is $2^s\le 2^{\omega(N)}$,
with $s$ the number of $p_i|N$ that split in $\Zee[\omega]$.

Step 4 computes the possible generators of the primitive parts of ideals
of norm~$N$ in $\Zee[\omega]$. 
It is not completely optimized as it does not take into account that 
different roots of $f$ modulo $N_1$ may coincide modulo $N_0$, and give 
rise to the same ideal $(N_0, \omega-r)$ in Step 4a.
It also unnecessarily treats the complex conjugate $(N_0, \omega-r')$
of every ideal $(N_0, \omega-r)$, whose generators (if any) are of
course the complex conjugates of the generators of $(N_0, \omega-r)$.

\head 4. Heuristic run time analysis
\endhead

\noindent
In this section, we present a heuristic run time analysis of
the Algorithm in the previous section, and numerical data
supporting this analysis.
\smallskip\noindent
{\it Assumption 1.}
For the elements $\alpha=k\alpha_0\in\Zee[\omega]$ of norm $N$ that we find in
Step~4a of our Algorithm, the norm of $1-\alpha$ will be an element of the 
Hasse interval ${\cal H}_N$ that, apart from being congruent to
$1\mod k$, does not appear to have any predictable primality properties.
Based on the prime number theorem, a reasonable assumption
is therefore that for varying $d$, $r$ and $N_0$, the norms found 
in Step 4b will be prime with `probability' at least $1/\log N$.
In other words, the number of times we expect to execute
Step 4b of our Algorithm before we find a prime value
is of order of magnitude $\log N$.
\smallskip\noindent
{\it Assumption 2.}
The input for Step 4b is provided by Step 4a, which finds the
generators of those ideals of norm $N$ in $\Zee[\omega]$ that
are principal.
The likelihood for a `random' ideal in $\Zee[\omega]$ to be principal
is $1/h_d$, with $h_d$ the class number of the ring of integers
$\Zee[\omega]\subset \Que(\sqrt{-d})$.
As we have no indication that the primitive ideals of norm $N_0$ arising
in Step 4a behave differently from random ideals in $\Zee[\omega]$, 
it seems reasonable that they will be principal with `probability'
around $1/h_d$.
\smallskip\noindent
The class number $h_d$ behaves somewhat irregularly as a function of $d$,
but its growth rate $d^{{1\over 2}+o(1)}$ was already found by Siegel.
In order to bound the number of times we execute the steps 4a and 4b, 
we need to bound the integers $d$ we encounter in Step 2, i.e.,
to find an upper bound $B_N$ for the minimal integer $d$ that occurs
in a solution to Problem 3.
Clearly, such an upper bound will be of heuristical nature, based on
the two `randomness assumptions' above.
As our Algorithm consists of a loop over $d=1, 2, 3, \ldots$,
and $d$ has to be factored in Step 2 to find if it is squarefree,
the value of $B_N$ is of great importance in estimating the run time,
and the success of our method depends on $B_N$ being `small' as a function
of $N$.
\medskip\noindent
{\bf Elliptic curves of prime order.}
In the case our input number $N$ is prime,
our Algorithm is similar to the first step of the
elliptic curve primality proving algorithm ECPP.
On input $N$, this algorithm looks for an imaginary quadratic field $K$
of small discriminant containing an element $\alpha$ of norm $N$
with the property that $N_{K/\Que}(1-\alpha)=N+1-\Tr_{K/\Que}(\alpha)$ is
{\it twice} a probable prime number $N'$.
If $\alpha\in K$ is found, $N$ becomes the order of the finite field $\FF$
and $2N'$ the number of points of an elliptic curve over $\FF$.
As $\#\FF$ and $\#E(\FF)$ occur symmetrically in all considerations,
this problem is almost identical to our Problem 3.
In fact, since finding a prime around a large number $N$ is heuristically
just as difficult as finding twice a prime around~$N$, the heuristic run
time for our Algorithm on prime input $N$ is {\it identical\/} to the
heuristic run time for the first step of ECPP on input $N$.
In accordance with the results in [\Mo, Section 3], we obtain 
the following.
\proclaim
{4.1. Theorem}
Let $N$ be a prime number.
Under the heuristic Assumptions 1 and~2,
the integer $d$ solving Problem 3 is of size $\widetilde O((\log N)^2)$, and
our Algorithm can be expected to find it in time $O((\log N)^{4+\varepsilon})$.
\endproclaim
\proclaim
{4.2. Corollary}
Under the heuristic Assumptions 1 and 2, Problem 2 admits a solution in time
$O((\log N)^{4+\varepsilon})$ for prime values of $N$.
\endproclaim
\noindent{\bf Proof of 4.2.}
We first use our Algorithm to find $d$, $\alpha$ and $p=N-1+\Tr(\alpha)$
solving Problem 3 for $N$;
the time $O((\log N)^{4+\varepsilon})$ needed for this dominates
the steps that follow.
We then construct the class polynomial $P_D$ for $D=\text{disc}(\Que(\sqrt{-d}))$
in time $\widetilde O(d)=\widetilde O((\log N)^2)$.
As $P_D$ has degree $h_d\approx\sqrt d$, finding a root $j$ of $P_d$ in
$\Fp$ takes time $\widetilde O(\deg(P_d) (\log p)^2)=\widetilde O((\log N)^3)$
[\MCA, Section 14.5].
An elliptic curve $E$ with $j$-invariant $j$ and its quadratic twist $E'$
will have $N=p+1-\Tr(\alpha)$ or $p+1+\Tr(\alpha)$ points.
Matching the group order with the curve can be done efficiently by
determining which of the two quantities annihilates random points
on the curve.
We know that only one of them does for either $E$ or $E'$
for all $p>229$ by [\SC, Theorem 3.2].
\hfill$\square$
\medskip\noindent
\noindent{\bf Proof of 4.1.}
For prime input $N$, our algorithm is rather simple.
For increasing values of~$d$,
it singles out those $d$ for which $N$ is not inert
in $\Zee[\omega_d]$ in Step 3; in Step 4,
it computes the primes over $N$ in $\Zee[\omega_d]$ and determines 
whether these are principal with a generator $\alpha$
for which $1-\alpha$ is a prime element.

The ring $\Zee[\omega_d]$ contains elements $\alpha$
of norm $N$ if and only if $N$ splits into principal primes of norm $N$.
For primes $N$ coprime to $2d$, this means that $N$ has to split 
completely in the Hilbert class field $H_d$ of $\Que(\sqrt{-d})$.
Our Assumption~2, which states that primitive ideals of norm $N$ should be
principal in $\Zee[\omega]$ with `probability' $1/h_d$, now reminds us of
the Chebotarev density {\it theorem}, which tells us that one out of every
$2h_d=[H_d:\Que]$ primes splits completely in $H_d$.
For $d>3$, it leads us to expect with `probability'
$1/(2h_d)$ that there are (up to conjugation)
exactly two integral elements $\alpha$ and $-\alpha$ of norm $N$.
With complementary probability $1-(2h_d)^{-1}$, there are no
elements of norm~$N$.
Thus, a value $d$ can be expected to yield an `on average' number
of $1/h_d$ elements of norm~$N$.

The average statement that the number of algebraic integers $\alpha\in
\Que(\sqrt {-d})$ of norm $N$ is asymptotically a fraction $1/h_d$
of the pairs $(d, N)$ tried is implied by Chebotarev's theorem in case
we fix $d$ and let the prime $N$ vary.
We are however in the case where $N$ is fixed and $d$ varies.
This is certainly different, but for varying $d$ up to a bound $B$
that is {\it small\/} with respect to $N$, it is Assumption 2
that we will find approximately $\sum_{d<B} 1/h_d$ elements of fixed norm $N$.
This is reasonable, provided that the fields $H_d$ are `close' to being
linearly independent over $\Que$.

It is not exactly true that the Hilbert class fields $H_d$ for the squarefree
integers $d<B$ we encounter form a linearly disjoint family of number
fields: the genus fields $G_d\subset H_d$ have many non-trivial
intersections.
However, in this family of fields, which has about $(6/\pi^2)B$ elements, 
there is a subfamily of fields $H_d$ coming from the prime numbers
$d\congr3\mod 4$ that is linearly disjoint over $\Que$.
This follows from the fact that for these primes $d$, the field $H_d$
is ramified only at $d$, so every field $H_d$ is linearly disjoint from
the compositum of the other fields $H_d$ in the subfamily.
As the given subfamily has asymptotically $B/(2\log B)$ 
elements, we can treat the family of fields $H_d$ with $d<B$ as being
linearly independent at the cost of allowing for lower order 
(logarithmic) factors in our estimates.
We can estimate the asymptotic size of the sum
$\sum_{d<B} 1/h_d$ for squarefree $d<B$ to be a positive constant times
$
\sum_{0<d<B} {1\over\sqrt d}\approx \int_0^B {dt\over \sqrt t}= 2\sqrt B.
$

We find that for $B$ tending to infinity, Assumption 2 implies that
the number of elements of prime norm $N$ coming from $d<B$
is bounded from below by some universal constant times $\sqrt B/\log B$.
By Assumption 1, we expect to need about $\log N$ elements of norm $N$ in Step 5b.
Thus, for prime values $N$ tending to infinity,
the size $B_N$ of the minimal $d$ solving Problem 3 can be expected to
be of size $\widetilde O((\log N)^2)$.
Note that $B_N$ is small with respect to $N$, as required in our
heuristical argument.

For the run time of the algorithm, we obtain $O((\log N)^{4+\varepsilon})$
exactly as in~[\Mo].
The main term in the run time comes from computing
$\widetilde O((\log N)^2)$ values of $\sqrt{-d}\mod N$, which each take
time $O((\log N)^{2+\varepsilon})$,
and from proving (as in [\Bern]) that the output is correct, i.e., that
we have found $\alpha$ of norm $N$ for which $N+1-\Tr(\alpha)$ is
indeed prime.
\hfill$\square$
\medskip\noindent
{\sl Numerical support.}
The table below shows the number of solutions 
$x,y\in\Zee_{\geq 1}$ to the equation $x^2+dy^2=4N$ for $d$ ranging over
the squarefree integers $d\in [1,B]$ for various $B$.
For $N$ we took the 5 primes following $10^{100}$ and $10^{200}$. 
Note that the spacing of primes around $10^{100}$ and $10^{200}$
is in accordance with Assumption 1.
$$
\vbox{\offinterlineskip
\halign{&\vrule#&\strut\quad#\quad\hfil\cr
width0pt &$\downarrow N$\hfill$B\to$&\hskip1pt\vrule&\hfil$1000$&&\hfil$4000$&&\hfil$16000$&&\hfil$64000$\cr
height2pt width0pt&\omit&\hskip1pt\vrule&\omit&&\omit&&\omit&&\omit\cr
\noalign{\hrule}
height1pt width0pt&\omit&\hskip1pt\vrule&\omit&&\omit&&\omit&&\omit\cr
\noalign{\hrule}
height2pt width0pt&\omit&\hskip1pt\vrule&\omit&&\omit&&\omit&&\omit\cr
width0pt &$p_1=10^{100}+267$&\hskip1pt\vrule&\hfil30&&\hfil57&&\hfil125&&\hfil232\cr
width0pt &$p_2=10^{100}+949$&\hskip1pt\vrule&\hfil41&&\hfil87&&\hfil161&&\hfil304\cr
width0pt &$p_3=10^{100}+1243$&\hskip1pt\vrule&\hfil22&&\hfil51&&\hfil93&&\hfil173\cr
width0pt &$p_4=10^{100}+1293$&\hskip1pt\vrule&\hfil39&&\hfil72&&\hfil145&&\hfil316\cr
width0pt &$p_5=10^{100}+1983$&\hskip1pt\vrule&\hfil29&&\hfil57&&\hfil123&&\hfil245\cr
width0pt &$q_1=10^{200}+357$&\hskip1pt\vrule&\hfil46&&\hfil91&&\hfil190&&\hfil354\cr
width0pt &$q_2=10^{200}+627$&\hskip1pt\vrule&\hfil24&&\hfil51&&\hfil98&&\hfil210\cr
width0pt &$q_3=10^{200}+799$&\hskip1pt\vrule&\hfil24&&\hfil47&&\hfil90&&\hfil184\cr
width0pt &$q_4=10^{200}+1849$&\hskip1pt\vrule&\hfil47&&\hfil81&&\hfil170&&\hfil376\cr
width0pt &$q_5=10^{200}+2569$&\hskip1pt\vrule&\hfil73&&\hfil140&&\hfil275&&\hfil532\cr}
}
$$
We see that the growth rate is indeed roughly proportional to $c_N\sqrt{B}$,
for some constant $c_N$: the numbers double if we quadruple $B$.

The data show that the size of $N$, when large with respect to $B$,
is irrelevant: only the class of the primes over $N$ in the class group
of $\Zee[\omega]$ is important, not the size of $N$.

Figure 1 below shows the number of solutions for $p_2$ and $p_3$. Inspecting
the data, we see that the growth rate is indeed close to $\sqrt{B}$. The
fluctuation in the graphs is caused by the somewhat irregular behaviour
of $h_d$. On a logarithmic scale, the graphs do look like straight lines
with slope $1/2$, see Figure 2.
\medskip\noindent
\centerline{\epsffile{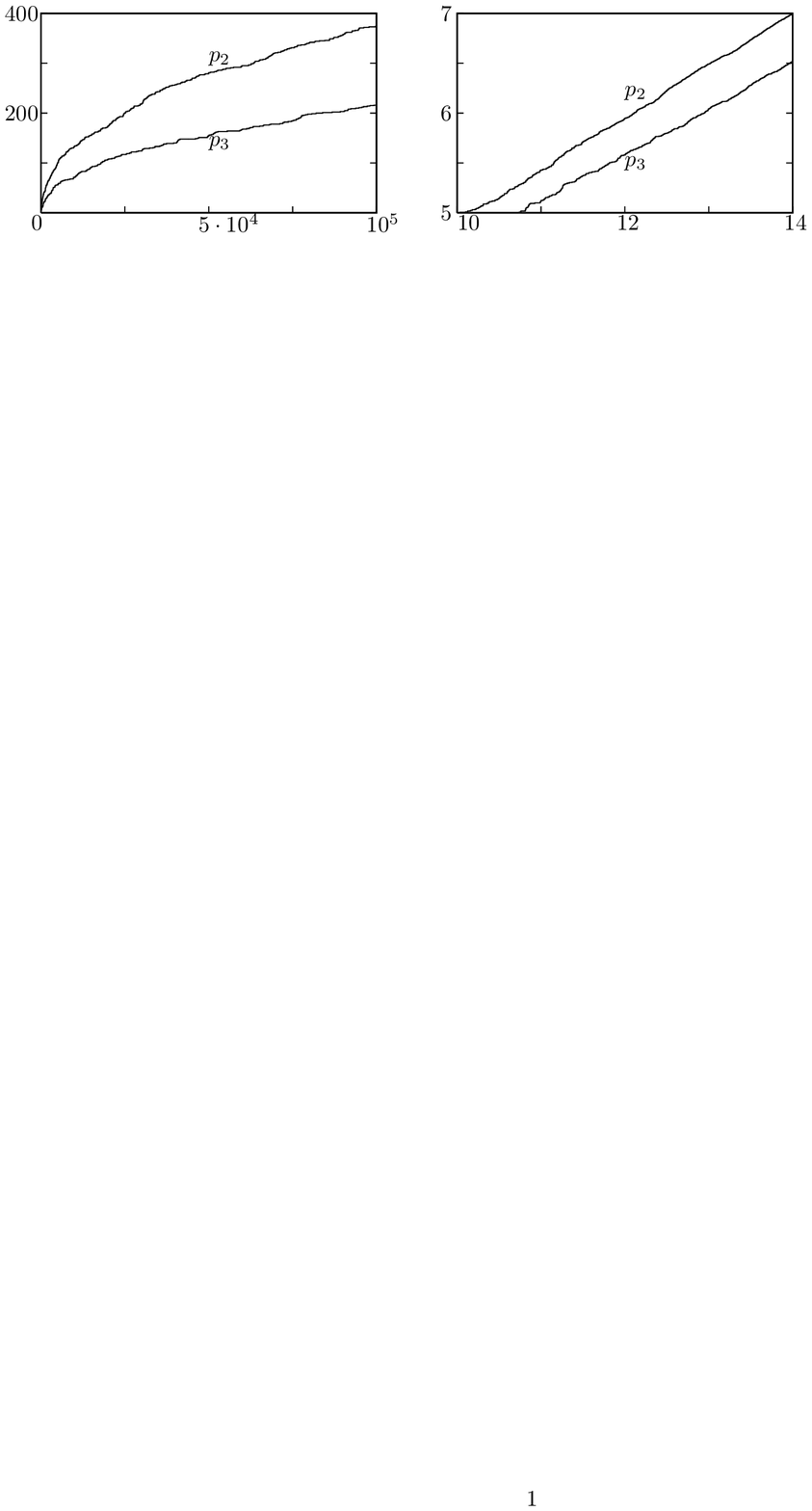}}
\par
\hskip2cm{\bf Figure 1}\hskip5cm{\bf Figure 2}\par
\medskip\noindent
There are clear differences in the constants $c_N$ for various $N$.
These can be explained by looking at the
contributions coming from composite $d$, which we could afford
to neglect in our analysis, but which play an important role
in practical situations.
For solvability of (2.4), it is clear that $N$ has to be a square
modulo all primes dividing~$d$. For even $d$, we also have the
condition ${N \overwithdelims() 2}=1$.
If we have ${N\overwithdelims() p}=1$ for many small primes $p$, there will 
most likely be more composite $d$ yielding solutions to (2.4). 
The most striking difference in the table occurs for $p_3$ and $q_5$. 
Looking at the Kronecker symbols ${q_5 \overwithdelims() p}$ for the first
eight primes $p \leq 20$, we only have ${q_5 \overwithdelims() p} = -1$ for
$p=3,11$. For $p_3$ this occurs for $p=2,3,5,13,17$. This
explains why $q_5$ `outperforms' $p_3$.  
The differences in the constants $c_N$ disappear if we only consider
primes $d \equiv 3 \bmod 4$ in our table. For $p_3$ we get $53$ solutions
up to $B=64000$ in this case and for $q_5$ we get $50$ solutions.

Whereas the number of generators of norm $N$ found in Step 5a for $d<B$
increases regularly, and roughly proportional to $\sqrt B$, 
Assumption~1 tells us that the number of times we have to test for primality in
Step 5b before we hit a prime number is $\log N$ {\it on average}.
As a consequence, we expect that the minimal $d=d(N)$ solving Problem 3 
is of size $O((\log N)^{2+\varepsilon})$, but not that $d(N)$
increases very regularly with~$N$ for prime values $N$.
For instance, the primes $p_1$ and $p_5$ above have rather similar
curves exhibiting the number of solutions found in Step 5a,
but the corresponding minimal discriminants 643 and 303267 are
quite far apart: they are the smallest and largest values found for the $p_i$.
However, the average value of~$d$ for the first 100 primes larger
than $10^{100}$ and the first 100 primes larger than $10^{200}$ are 
$82170\approx(\log(10^{100}))^{2.08}$ and $396030\approx(\log(10^{200}))^{2.10}$,
respectively.
Their quotient 4.8 is not too far from the factor $4$ we expect.
\medskip\noindent
{\bf Elliptic curves of arbitrary order.}
The Assumptions 1 and 2 at the beginning of the section
also provide a heuristic run time analysis for arbitrary input $N$.

Assume first that $N$ is squarefree, say
$N=\prod_{i=1}^{\omega(N)} p_i$ with $p_i$ prime.
In Step 3a, all $d$ are discarded for which one of the primes $p_i$
is inert in $\Zee[\omega_d]$, so we will only be working in Step 4 with
those $d$ for which none of the $\omega(N)$ Kronecker symbols
$D\overwithdelims()p_i$ equals $-1$.
This can be a set of integers of density as small as $2^{-\omega(N)}$
inside the set of all squarefree integers,
and in case $N$ is in the zero-density subset of integers satisfying
the equivalent inequalities
$$
2^{\omega(N)}> (\log N)^2
\Longleftrightarrow
\omega(N)> {2\over \log 2} \log\log N = 2.88539 \log\log N
$$
it is clear that we can no longer expect the integer $d$
solving Problem 3 to be of size at most $(\log N)^{2+\varepsilon}$.

Despite the scarcity of suitable $d$ for large values of $\omega(N)$,
it is still the case that we expect the number of elements of norm
$N$ coming from $d<B$ to grow at least as fast as some
universal constant times $\sqrt B/\log B$ if $B$ tends to infinity.
Indeed, looking as before at the prime numbers $d\congr3\mod 4$ (not dividing
$N$) up to $B$, we see that there are ideals of norm $N$ only for a fraction
$2^{-\omega(N)}$ of them.
However, for each $d$ meeting the $\omega(N)$ quadratic conditions,
the number of ideals $I$ of norm~$N$ equals $2^{\omega(N)}$: we can take
$I=\prod_{i=1}^{\omega(N)} \gothp_i$, with $\gothp_i$ one of the two
primes dividing~$p_i$ in~$\Zee[\omega_d]$.
This means that the growth with $B$ of the number of ideals of
norm $N$ coming from $d<B$ is {\it independent\/} of the value of
$\omega(N)$: with increasing $\omega(N)$ they occur for fewer $d$,
but the decrease in contributing $d$ is exactly compensated by the
number of ideals provided by such $d$.
Our expected number of elements of norm
$N$ coming from $d<B$ is therefore unchanged with respect to the
case of primes $N$ discussed before.

The problem with the asymptotic growth $\sqrt B/\log B$ of elements
of norm $N$ coming from a thin subset of $d<B$ is that $B$ may
have to be large to observe this growth rate: clearly the expected number
$2^{-\omega(N)}B$ of contributing $d<B$ should not be too small.
As we want to take $B\approx (\log N)^2$, we can only use our previous
estimate for the expected size of the integer $d$ solving Problem 3
in the case $2^{\omega(N)}\ll (\log N)^2$.
In the `opposite' case $2^{\omega(N)}\gg (\log N)^2$, finding a {\it single\/}
quadratic ring $\Zee[\omega_d]$ in which all primes $p_i|N$ split
completely is what the Algorithm needs to achieve: there will be
$2^{\omega(N)}$ ideals of norm $N$ in this ring, of which Assumption~2
tells us we can expect $2^{\omega(N)}/h_d\approx 2^{\omega(N)}/\sqrt d$
to be principal.
As the smallest~$d$ satisfying the $\omega(N)$ quadratic conditions
imposed by the $p_i$ is expected to be of order of magnitude $2^{\omega(N)}$,
we will find $2^{\omega(N)/2}\gg \log N$ elements $\alpha$ of norm $N$
in $\Zee[\omega_d]$.
By Assumption 1 this will lead to a prime element $1-\alpha$.
\proclaim
{4.3. Theorem} 
Under the heuristic Assumptions 1 and~2,
the integer $d$ solving Problem 3 is of size
$\widetilde O ((\log N)^2 + 2^{\omega(N)})$, and our 
Algorithm can be expected to find it in time 
$O(2^{\omega(N)}(\log N)^{4+\varepsilon})$.
\endproclaim
\proclaim 
{4.4. Corollary}
Under the heuristic Assumptions 1 and 2, Problem 2 admits a solution in time
$O(2^{\omega(N)}(\log N)^{4+\varepsilon})$.
\endproclaim
\noindent{\bf Proof of 4.4.}
Analogous to the proof of 4.2.\hfill$\square$
\medskip
\noindent{\bf Proof of 4.3.}
We saw that for squarefree $N$, the size of the integer $d$
solving Problem~3 is of size $\widetilde O ((\log N)^2)$ in case
$2^{\omega(N)}$ is smaller magnitude. If it is bigger, the term
$2^{\omega(N)}$ becomes dominant and determines the expected size
$\widetilde O (2^{\omega(N)})$ of $d$.

If $N$ is not squarefree, the Algorithm has an increased number of 
possibilities to find ideals and elements of norm $N$ for each value of $d$.
Primes occurring to even exponents are no longer an obstruction if they
are inert in $\Zee[\omega_d]$: they get absorbed in $k_1$ in Step 3 and
no longer occur in $N_1$ in Step 4.
Splitting primes occurring to higher exponents lead to square divisors
$k_2^2|N_1$ in Step 4, and to various ideals $(N_0, \omega-r)$ that can
be tested for principality in Step 4a.
The extra ways to find elements of norm $N$ is an advantage as it will lead 
to a smaller bound $B_N$ for the minimal $d$ solving Problem 3.
In particular, $B_N$ will be of size $\widetilde O ((\log N)^2 + 2^{\omega(N)})$
for {\it all\/} $N$.

In order to estimate the run time of the Algorithm, we observe that by
Assumption 1, Step 4b will be executed about $\log N$ times 
until a probable prime norm is found, and a true primality proof
taking expected time $O((\log N)^{4+\varepsilon})$ is needed.
This is the dominant term in the time spent on Step 4b.
The number of times Cornacchia's algorithm is executed in Step 4a
to yield the $\log N$ generators going into Step 4b is by Assumption 2
no more than $O(\sqrt {B_N} \log  N)$, as the class numbers $h_d$ for $d<B_N$ are 
no bigger than $\sqrt {B_N}$. 
As Cornacchia's algorithm takes time $O((\log N)^{2+\varepsilon})$,
we expect to spend time  $O(\sqrt {B_N}(\log N)^{3+\varepsilon})$ in Step 4a.

In order to find the roots $(r\mod N_1)$ of $f$ in Step 4, we first extract
the square roots $\sqrt{-d}$ modulo each of the primes $p_i$ that split
in $\Zee[\omega_d]$, in time at most $O(\omega(N) (\log N)^{2+\varepsilon})$.
For each choice of square roots, there is a root $(r\mod N_1)$ of $f$ that
can be found using the Chinese remainder theorem, in time 
$\omega(N)(\log N)^2$.
Each time we apply the Chinese remainder theorem, we use the root 
$(r\mod N_1)$ obtained in Cornacchia's algorithm in Step 4a.
The number of times we apply the Chinese remainder theorem is therefore
bounded by the number of times $O(\sqrt {B_N} \log  N)$ we apply
Cornacchia's algorithm. 
We find that the total time spent on finding roots
$(r\mod N_1)$ is no more than 
$O(\sqrt {B_N} \omega(N) (\log N)^3)$.
Taking all parts of Step 4 together, the total time spent in
Step 4 becomes $O(\sqrt {B_N} \omega(N) (\log N)^{3+\varepsilon})$.
This is $O((\log N)^{4+\varepsilon})$ in the case
$2^{\omega(N)}\ll (\log N)^2$, and $O(2^{\omega(N)/2}(\log N)^{4+\varepsilon})$
in general.

Outside Step 4, no substantial computing is done, only some administration for
the relatively small integer $d$, which takes values up to ${B_N}$.
In cases where $B_N$ is of order of magnitude $2^{\omega(N)}\gg (\log N)^2$,
doing this administration is not negligeable because of
the large number of values taken on by $d$. 
Taking this into account, we find that the heuristic run time is
bounded in all cases by $O(2^{\omega(N)}(\log N)^{4+\varepsilon})$.
\hfill$\square$
\medskip\noindent
{\sl Numerical support.}
Figure 3 below shows how the number of solutions $x,y \in \Zee_{\geq 1}$
to the equation $x^2+dy^2=4N$ for $d$ ranging over
all squarefree integers $d\in [1,B]$ varies with $B$
for different number $\omega(N)$ of prime factors of $N$.
The graphs are given for $N=N_1,N_2,N_3,N_{10}$, where $N_k$ is the product of
the first $k$ primes larger than $10^{10}$. 

We see that the graphs for $N_1,N_2$ and $N_3$ behave quite similarly. This
is what we expected if the number of solutions is independent of $\omega(N)$.
The graph for $N_{10}$ appears to be quite different from the others,
and this is because $2^{\omega(N_{10})}=2^{10} = 1024$ is here of
the same order of magnitude as the values of $B$ in the graph.
There are here fewer $d$ for
which we have a solution to $x^2+dy^2=4N_{10}$, but if we do have a solution,
we immediately get many. For instance, the first `jump' in the graph 
occurs for the prime value $d=1949$ and we get 28 solutions for this $d$.
This is in nice accordance with the heuristics, which tell us to
expect the first solutions to occur for around $d \approx 2^{10} = 1024$,
and to be about $2^5 = 32$ in number.
\medskip\noindent
\centerline{\epsffile{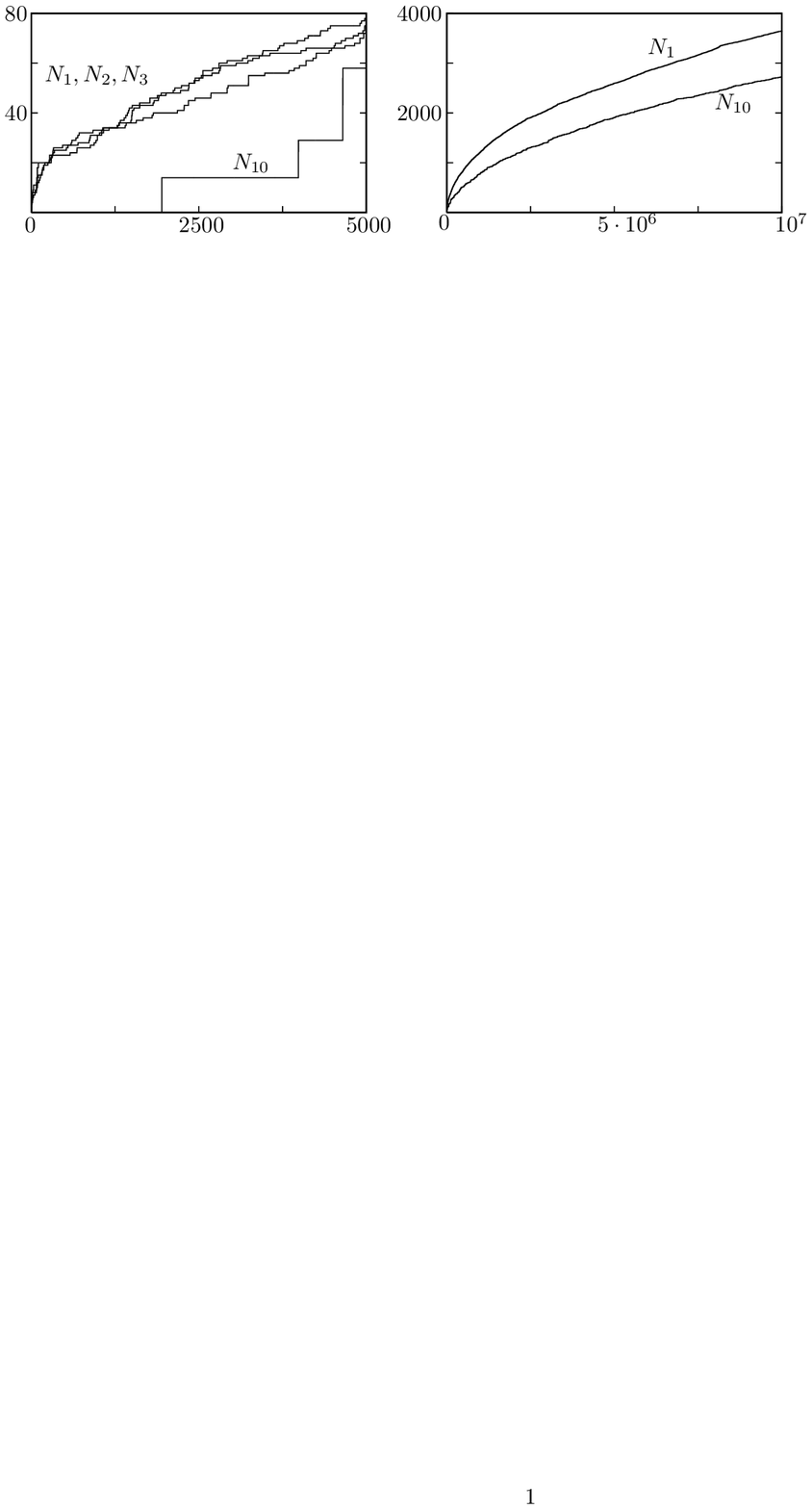}}
\par
\hskip2cm{\bf Figure 3}\hskip5cm{\bf Figure 4}\par
\medskip\noindent
The irregularity of the graph for $N_{10}$ disappears if we look at
values of $B$ that are large in comparison to $2^{\omega(N_{10})}$.
Figure 4 shows the graph for $N_{10}$ for $B$ up to $10^7$. 
It is now similar in nature to that of $N_1$, and exhibits the familiar
$\sqrt{B}$-profile.

The graph in Figure 5 below illustrates the dependence on the number of square 
divisors of $N$. It shows the number of solutions for $N_1$,
$3^2\cdot N_1$, $3^2\cdot 5^2\cdot N_1$ and $3^2 \cdot 5^2 \cdot 7^2\cdot N_1$. 
If $N$ has square divisors, we potentially test the principality of 
more ideals in step 4 of our Algorithm, so we expect to obtain more solutions.
Replacing $N_1$ for example by $3^2\cdot N_1$, we expect to get on average
a double amount of solutions for $d\congr 1\mod 3$.
The gain is a constant factor $>1$ that increases with the amount of
square divisors of $N$.
\medskip
\centerline{
\epsffile{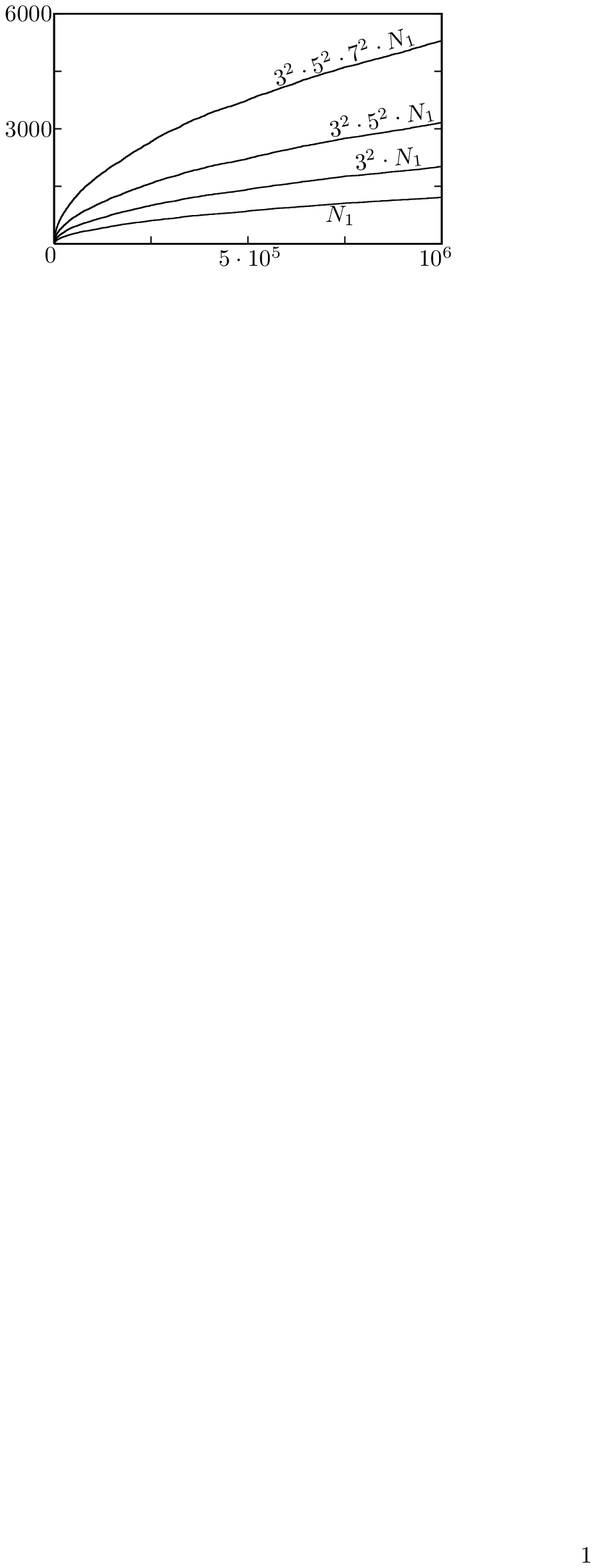}}
\centerline{\bf Figure 5}

\head 5. Examples and practical considerations
\endhead

\noindent
The description of the Algorithm in Section 3
is intended to facilitate the run time estimate in Section 4,
it does not address practical issues that are important
in computing large examples.
In this section, we explain how we find
solutions to Problem 2 form large values of $N$ that
are either prime or equal to a power of 10.
\medskip\noindent
{\bf Elliptic curves of large prime order.}
{}From the description of the algorithm we gave in the previous section, and more
in particular its relation to ECPP, it is clear that one should be able
to construct a curve having a large prime number $N$ of points in all
cases where ECPP, as described in [\Mo], can prove primality of a number
of the same size.
To do so, it makes sense to apply an idea attributed to J.~Shallit in [\Mo]
to speed up the computation.
This idea starts from the observation that for
large prime numbers $N$, the Algorithm spends a lot of time
in evaluating $(\sqrt{-d}\mod N)$ for all squarefree $d$ up
$B_N\approx(\log N)^2$ having ${-d\overwithdelims()N}=1$.
We noticed already in the previous section that if the equation
$$
x^2+dy^2=4N 
$$
admits integral solutions, then $N$ is a square modulo all primes
dividing $D=\text{disc}(\Que(\sqrt{-d})$. 
It reflects the fact that if $N$  splits completely
in the Hilbert class field $H_d$ of $K=\Que(\sqrt{-d})$, then 
it certainly splits completely in the genus field $G_d\subset H_d$ of $K$.
As $G_d$ is obtained by adjoining to $K$ the square roots of
$p^* = (-1)^{(p-1)/2}p$ for all odd prime divisors $p|d$,
we have ${p^*\overwithdelims() N}={N\overwithdelims() p}=1$ in this case.

Once we know that those $d$ providing solutions are essentially
products of primes having the right quadratic character with respect to
$N$, the idea suggests itself to look at those $d$ only that are
constructed as products of such primes.
Creating $d$ from a `basis' of primes
$p$ with $\bigl( {p^*\over N} \bigr)=1$ allows us to compute 
$\sqrt{p^*} \bmod N$ for such $p$, and store the values in a list.
For $p=2$, one uses the square roots of $-1$, $2$ and $-2$ that can
be extracted modulo $N$.
For each $d$ constructed from our basis of primes,
$\sqrt{-d} \bmod N$ can be obtained by multiplying the square roots
of primes modulo $N$ we stored.
Considering only products of two primes from our basis allows
us to reduce the number of square root extractions modulo $N$ from
$O((\log N)^2)$ to $O(\log N)$, at the expense of extra multiplications
modulo $N$ and an increased storage requirement.
In practice, we consider $d$ with at most $3$ prime divisors.
One thing we lose in this approach is the guarantee that we really find the
{\it smallest\/} solution $d$ to Problem 3.
\medskip\noindent
{\bf 5.1. Example.}
Take $N=\text{nextprime}(10^{2004}) = 10^{2004}+4863$, the exponent
2004 being the year we found our method.
For this $N$, we have $\log(N)=4614.3$ and $(\log(N))^2=2.13\cdot 10^7$.
There are 324 primes $p$ less than $5000$ with $\bigl( {p^* \over
N} \bigr)=1$, and we compute and store $\sqrt{-1} \bmod N$ 
and all square roots $\sqrt{p^*} \bmod N$.
We now have ${325\choose 3} =5668650$ squarefree values of $d$ at our
disposal having up to 3 prime divisors from our base, and we know $N$
to split completely in all genus fields $G_d$. 

The $104415$-th value of $d$ we tried was $d = 59\cdot 523 \cdot 2579=79580203$.
For this value of $d$, we found a solution
$$
x=1885782\ldots693127
$$
to $x^2+dy^2=4N$ for which
$$
p = N+1-x=999999\ldots99999811421\ldots8311737
$$
is a 2004-digit prime. 
In each case, the dots represent 990 digits that we omitted.

The class polynomial $P_{-d}$ has degree $1536$ and coefficients up to 41984 
digits. Modulo $p$, the polynomial $P_{-d}$ splits completely.
Taking $j$ to be the smallest positive integer satisfying 
$P_{-d}(j)\congr 0\mod p$
we put $a = {27j \over 4(1728-j)}\in\Fp$.
Then the curve given by
$$
E_a :\thinspace\thinspace  Y^2 = X^3 + aX -a
$$
has CM by $\O_{-d}$. As the point $(1,1)\in E_a(\Fp)$ 
does not have order $N$, the quadratic twist 
$E_a': Y^2 = X^3 +9aX-27a$ of $E_a$ has $N$ points.
This can be verified by picking a random point $P\in E_a'(\Fp)$
and checking that we have $N\cdot P=0$.

The value of $d$ we find here is in fact the smallest $d$ solving
Problem 3 for our $N$.
Our algorithm did 565 primality tests before we found the
solution above.
Finding $d$ and $p$ took about 10 minutes on our standard PC,
and another 3 hours were needed to find and factor $P_{-d}$.
Once we find $j$, the final result is almost immediate.
If we trust the input value $N$ as being a true prime number, 
there is no need to prove that $p$ is prime.
As in ECPP, this follows from the fact that $E'$ has a non-trivial
point that is annihilated by $N$.
\medskip\noindent
{\bf Elliptic curves of 10-power order.}
We indicated in our analysis in Section 4 that for input values of $N$ having
a large number of square divisors, the integer $d$ solving Problem 3
will be much smaller than the upper bound for squarefree $N$
occurring in Theorem 4.3.
This can be illustrated by looking at the values $N=10^k$ for $k\ge1$,
which have $\log N\approx 2.3 k$.
As none of the prime divisors 2 and 5 of $N$ is inert in the field
$\Que(i)$ and the prime 5 is split, there are already
many solutions to the norm equation $x^2+y^2=N$ for the very first value $d=1$.
In fact, as we have $h_d=1$ there is no need for a Cornacchia algorithm,
and the elements of norm $N=2^k5^k$ in $\Zee[i]$ are the $4k+4$ elements
$\alpha_{s,t}= i^s (1+i)^k (2+i)^t (2-i)^{k-t}$ with $s\in\{0, 1, 2, 3\}$
and $t\in\{0, 1,\ldots, k\}$.
Up to conjugacy, we have about $2k=.87 \log N$ elements, so we expect that
for a positive fraction of all $k$-values,
$d=1$ gives~rise to a prime $p$ and a twist $E$ of the curve $Y^2=X^3+X$
having exactly $10^k$ points over $\Fp$.
As the graph below indicates, this fraction appears to be
close to 0.92.

\centerline{
\epsffile{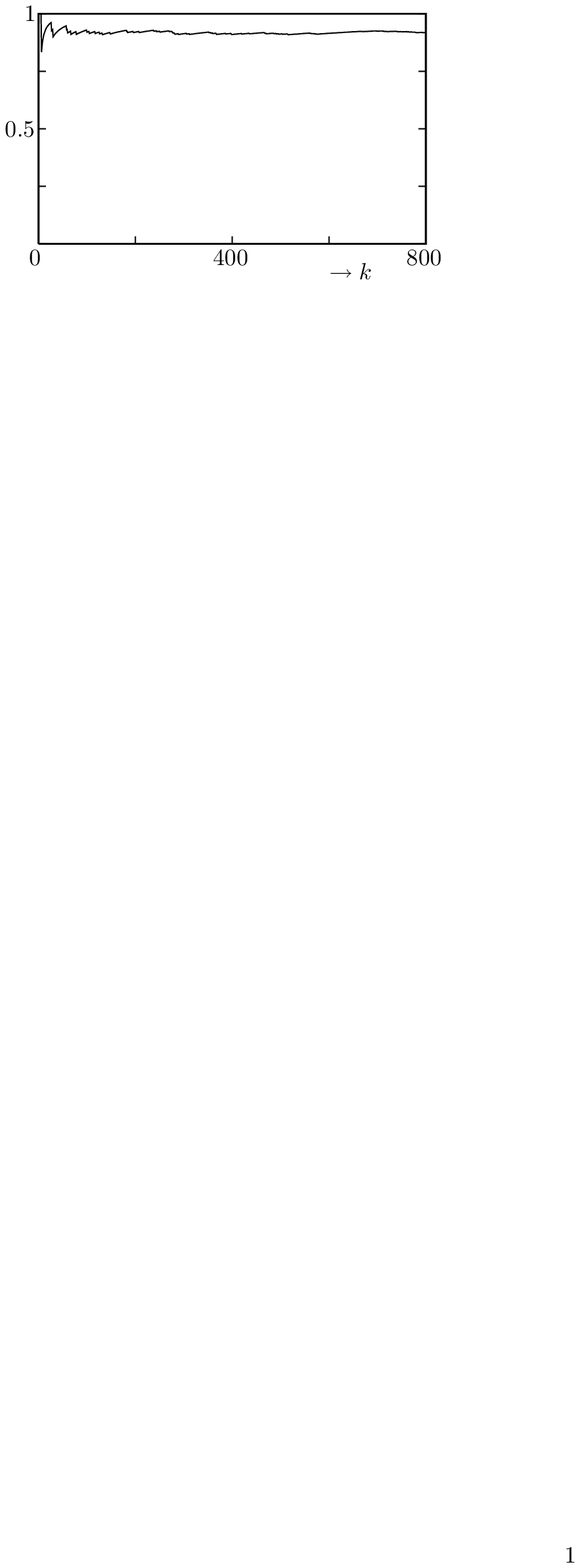}}
\centerline{\bf Figure 6}
\medskip\noindent
{\bf 5.2. Example.}
Take $k=2004$. We find that for $(s,t)=(2,499), (0,527), (0,671)$,
the element $\alpha_{s,t}=i^s (1+i)^{2004} (2+i)^t (2-i)^{2004-t}$
of norm $10^{2004}$ has the property
that $p=N_{\Que(i)/\Que}(1-\alpha_{s,t})$ is prime. 
The curve $Y^2=X^3+X$ having $j=0$ and CM by $\Zee[i]$ has
4 twists over $\Fp$ for each of these $p$, but in all cases
$Y^2=X^3+X$ is the curve having $10^{2004}$ points.
This follows from a result in [\CF] going back to Gauss.
It says if we choose the prime element $\pi=a+bi$ dividing a prime
$p\congr 1\mod 4$ in $\Zee[i]$ to satisfy $\pi \equiv 1 \bmod (1+i)^3$,
then the curve $Y^2=X^3+X$ has exactly 
$
p+1- {-1\overwithdelims() \bar\pi}_4 (\pi+ \bar\pi)= p+1-2i^{1-a}a
$ 
points over $\Fp$.
In our case, $\pi=1-\alpha_{s,t}$ and $a$ are congruent to 1 modulo 
$(1+i)^{2004}=-2^{1002}$, so we already know that $Y^2=X^3+X$ is
the right curve before actually computing $p$.
\medskip\noindent
For the purpose of constructing curves having $N=10^k$ points,
there are small values of $d$ that conjecturally work for almost all
values of $k$, not just for a positive fraction of them.
These $d$ have the property that 2 and 5 {\it both\/} split completely
in $\Que(\sqrt{-d})$, i.e., they satisfy $d\congr 31, 39\mod 40$.
For such $d$, the number of ideals of norm $N$  grows quadratically
in $k$, and hence in $\log N$. If we fix $d$, and hence $h_d$,
the number of elements of norm $N$ in $\Que(\sqrt{-d})$ will also
grow quadratically in $\log N$,
and our Assumption 2 implies that such $d$ will work for all
but finitely many $k$.
\medskip\noindent
{\bf 5.3. Example.}
Let $\rho$ be a zero of $X^3+X+1$. 
Then $\rho$ is the value of the Weber function $\goth{f}(z) = 
\zeta_{48}^{-1}\cdot {\eta({z+1\over 2})\over\eta(z)}$ at $-23-1/\omega_{31}$,
and a generator of the Hilbert class field of $\Que(\sqrt{-31})$.
An elliptic curve $E_j/\Que(\rho)$ having $j$-invariant
$j=(\rho^{24}-16)^3/\rho^{24}$ has endomorphism ring $\Zee[\omega_{31}]$.
We may take
$$
E_j: Y^2 = X^3 +3j(1728-j) + 2j(1728-j)^2
$$
which has good reduction outside $2,3,11,17,23,31$.
For all values $1 \leq k \leq 1000$ except $k=1,2$, there exist
primes of the form
$$
p=x^2+31y^2=10^k-1+2x.  \eqno{(5.4)}
$$
To find them, we write $(\omega_{31}+1)=\goth{p}_2\goth{p}_5$
and note that an $\Zee[\omega_{31}]$-ideal
$$
\goth{p}_2^s\cdot \bar\goth{p}_2^{k-s}\cdot \goth{p}_5^t\cdot \bar\goth{p}_5^{k-t}
$$
of norm $10^k$ is principal if and only if we have $s\equiv t \bmod 3$. We
use Cornacchia's algorithm to find the generators $\alpha$ for the principal 
ideals and test whether $N(1-\alpha)$ is prime.
For primes satisfying (5.4), either the reduction $\bar E_j/\FF_p$ of $E_j$
over a prime over $p$ in $\Que(\rho)$ or its quadratic twist
has exactly $10^k$ rational points over $\FF_p$.
It is likely that $k=1,2$ are the {\it only\/} values
of $k$ for which no prime $p$ of the form (5.4) exists,
but this is probably very hard to prove.

\enddocument